\documentclass[12pt]{amsart}
\usepackage{amscd}
\usepackage{verbatim}
\usepackage{amssymb, amsmath, amsthm, amscd,ifthen}
\usepackage[dvips]{graphics}
\usepackage[cp866]{inputenc}
\usepackage{graphicx}
\usepackage{epsfig}
\usepackage{color}

\usepackage{amsmath,amssymb,amscd,}
\usepackage[all]{xy}

\usepackage{amssymb}

\newtheorem{thm}{Theorem}[section]

\newtheorem{prop}[thm]{Proposition}

\newtheorem{Ex}[thm]{Example}

\newtheorem{lemma}[thm]{Lemma}
\newtheorem{cor}[thm]{Corollary}

\theoremstyle{definition}
\newtheorem{rem}[thm]{Remark}
\newtheorem{dfn}[thm]{Definition}

\newtheorem{prob}[thm]{Problem}

\title{Alexander $r$-tuples  and  Bier complexes }
\author[D. Joji\'{c}]{Du\v{s}ko Joji\'{c}}
\author[I. Nekrasov]{Ilya Nekrasov}
\author[G. Panina]{Gaiane Panina}
\author[R. \v{Z}ivaljevi\'{c}]{Rade \v{Z}ivaljevi\'{c}}

\address[D. Joji\'{c}]{ Faculty of Science, University of Banja Luka}
\address[I. Nekrasov]{{Chebyshev Laboratory}, St. Petersburg State University}
\address[G. Panina]{Mathematics \& Mechanics Department, St. Petersburg State University}
\address[R. \v{Z}ivaljevi\'{c}]{Mathematical Institute SASA, Belgrade}


\begin{document}

\begin{abstract}
{We introduce and study}  {\em Alexander $r$-tuples} $\mathcal{K}
= \langle K_i\rangle_{i=1}^r$ of simplicial complexes, as a common
generalization of { pairs of Alexander dual} complexes (Alexander
$2$-tuples) and $r$-unavoidable complexes of \cite{bfz}. In the
same vein, the {\em  Bier complexes}, defined as the deleted joins
$\mathcal{K}^\ast_\Delta$ of Alexander $r$-tuples, include both
standard {\em Bier spheres} and {\em optimal multiple chessboard
complexes} (Section~\ref{sec:multi-chess}) {as interesting,
special cases.}

Our main results are Theorem~\ref{ThmMain} saying that (1) the
$r$-fold deleted join of Alexander $r$-tuple is a pure complex
homotopy equivalent to a wedge of spheres, and (2) the $r$-fold
deleted join of a collective unavoidable $r$-tuple is
$(n-r-1)$-connected, and a classification theorem
(Theorem~\ref{thm:class} and Corollary~\ref{cor:Bier-cor}) for
Alexander $r$-tuples and Bier complexes.

\end{abstract}

\maketitle \setcounter{section}{0}

\section{Introduction}

Topological combinatorics utilizes  methods from algebraic (combinatorial)
 topology to solve problems in combinatorics and discrete geometry. Among
the highlights  (that strongly influenced the subsequent
developments),   and early achievements of topological combinatorics are
the solution of Kneser conjecture (L.~Lov\'{a}sz, 1978),
topological Tverberg theorem (I.~B\'{a}r\'{a}ny, S.B. Shlosman, A.~Sz\H ucs, 1981),
N.~Alon's `Splitting necklace theorem' (1987),
and many others, see \cite{Bjo95, Mat, Z04} for an overview and introduction.

\medskip
Simplicial complexes are among the central objects of study in topological combinatorics.
Their role in this subject can be compared to the role of manifolds in differential geometry and topology,
for illustration R.~Forman's `Discrete Morse theory' (Section~\ref{section_prelim})  exemplifies a fruitful interplay of ideas
and techniques from these areas.

\medskip

In this paper we introduce ``Alexander $r$-tuples of simplicial complexes'' and
closely related ``collective $r$-unavoidable complexes''
(Section~\ref{sec:collectively}), as unifying concepts that bring
together {\em Alexander pairs} of mutually dual complexes, and
{\em  $r$-unavoid\-able complexes} of Blagojevi\'{c}, Frick, and
Ziegler \cite[Definition~4.1]{bfz}.

\medskip
The deleted join operation, applied to an Alexander pair $(K,
K^\circ)$, yields a combinatorial sphere $Bier(K) = K\ast_\Delta
K^\circ$, known as the Bier sphere associated to $K$, see
\cite[Section~5.6]{Mat}. The special case of a self-dual complex
$K = K^\circ \subset 2^{[n]}$ is of particular importance. In this
case the Bier sphere $Bier(K) = K\ast_\Delta K$ is a
$\mathbb{Z}_2$-complex and its equivariant $\mathbb{Z}_2$-index is
${\rm Ind}_{\mathbb{Z}_2}(Bier(K)) = {\rm
Ind}_{\mathbb{Z}_2}(S^{n-2}) =  n-2$. This fact alone has many
interesting consequences, including the Van Kampen-Flores theorem
\cite[Theorem~5.1.1]{Mat} which says that the $d$-skeleton
$(\sigma^{2d+2})^{\leq d}$ of a $(2d+2)$-dimensional simplex is
non-embeddable in $\mathbb{R}^{2d}$.

\medskip
The {\em $r$-unavoidable complexes} \cite{bfz} play the central
role in applications of the `constraint method' of Blagojevi\'{c},
Frick, and Ziegler. This method, also known under the name
`Gromov-Blagojevi\'{c}-Frick-Ziegler reduction', has found
numerous applications to theorems of Tverberg-Van Kampen-Flores
type. We refer the reader to \cite[Section~2.9(c)]{Grom-2} and
\cite{bfz} for the original exposition of this beautiful technique
(see also our Section~\ref{sec:constraint} for a brief overview).

\medskip
The $2$-unavoidable complexes are easily identified as superdual
complexes $K\supseteq  K^\circ$. From here it easily follows that
self-dual complexes are precisely (inclusion) minimal
$2$-unavoidable complexes.

\medskip
Moreover, it was shown in \cite{jvz-3} (Theorem~3.6) that if $K$
is an $r$-unavoidable complex, then the associated $r$-fold
deleted join $K^{\ast r}_\Delta = K\ast_\Delta\dots\ast_\Delta K$
is a $S_r$-complex such that {the equivariant $G$-index} ${\rm Ind}_{G}(K^\ast_\Delta)\geq
n-r$ (where $r=p^k$ is a prime power and $G =
(\mathbb{Z}_r)^k\subset S_r$ is an elementary abelian group).

\medskip
The outline above leads to the conclusion that $r$-unavoidable
complexes can be interpreted as $r$-fold analogues (relatives) of
Alexander self-dual complexes, with many nice properties
preserved. It may be tempting to extend this analogy further, to
include $r$-fold generalization of (not necessarily symmetric)
Alexander dual pairs. The following research problem summarizes
the desirable properties of such an extension.

\begin{prob}\label{prob:main_problem}
Describe a property $\mathcal{P}_r$ of collections $\mathcal{K} =
\langle K_i\rangle_{i=1}^r = \langle K_1,\dots, K_r\rangle$ of
simplicial complexes on the same vertex set, $K_i\subset 2^{[n]}$,
such that:
\begin{enumerate}
 \item[(1)] If $r=2$ then a pair of complexes $\langle K_1,
 K_2\rangle$ satisfies $\mathcal{P}_2$ if and only if $\langle K_1,
 K_2\rangle$ is an Alexander superdual pair in the sense that $K_1\supseteq K_2^\circ$
 (equivalently $K_2\supseteq K_1^\circ$);

 \item[(2)] If $K_1=\dots= K_r=K$ then $\mathcal{K}$ satisfies
 $\mathcal{P}_r$ if and only if $K$ is an $r$-unavoidable complex;

 \item[(3)] If $\mathcal{K}\in \mathcal{P}_r$ then the deleted join
 $\mathcal{K}^\ast_\Delta = K_1\ast_\Delta\dots \ast_\Delta
 K_r$ is an $(n-r-1)$-connected complex.
\end{enumerate}
Moreover, it is desirable to describe a stronger property
$\mathcal{P}^\sharp_r\subset \mathcal{P}_r$ such that:

\begin{enumerate}
 \item[($1^\sharp$)]
$\langle K_1, K_2\rangle \in \mathcal{P}^\sharp_2$ if and only if
$K_1 = K_2^\circ$;
 \item[($2^\sharp$)] If $K_1=\dots= K_r=K$ and $\mathcal{K}\in \mathcal{P}_r^\sharp$,
 then $K$ is an (inclusion) minimal $r$-unavoidable complex;

 \item[($3^\sharp$)] If $\mathcal{K}\in \mathcal{P}_r^\sharp$ then the deleted join
 $\mathcal{K}^\ast_\Delta = K_1\ast_\Delta\dots \ast_\Delta
 K_r$ has the homotopy type of a wedge of $(n-r)$-dimensional spheres.
\end{enumerate}
\end{prob}

Motivated by Problem~\ref{prob:main_problem}, we describe
(Definition~\ref{dfn:pigeonhole}) the class $CU_r$ of ``collective
$r$-unavoidable complexes'', as our primary candidate for the
class $\mathcal{P}_r$. Individual $r$-unavoidable complexes often
arise from the `pigeonhole principle' (see \cite[Lemma~4.2]{bfz}).
For this reason we may occasionally say that an ordered collection
$\mathcal{K} = \langle K_i\rangle_{i=1}^r \in CU_r$ has the { {\em
pigeonhole property, }} or that $\mathcal{K}$ itself is a
\emph{pigeonhole $r$-tuple.}

\medskip

We introduce the class $\mathcal{A}_r$ of ``Alexander $r$-tuples
of simplicial complexes'' (Definition~\ref{DfnAlexander}), as the
most regular class of ``collective $r$-unavoidable complexes'',
and as our primary candidate for the class $\mathcal{P}^\sharp_r$.

\medskip
Finally,  {\em Bier complexes} (Section~\ref{sec:Bier-compl})
arise as the deleted joins of Alexander $r$-tuples, in perfect
analogy with the case of standard
 Bier spheres, which arise as deleted joins of Alexander
pairs of complexes.

\subsection{Summary of the main results}

The core of the paper are the results showing that the collective
$r$-unavoidable complexes (and their deleted joins) as well as the
Alexander $r$-tuples (and the associated  Bier  complexes) indeed
satisfy the properties listed in Problem~\ref{prob:main_problem}.
Perhaps the most interesting among them are the following {(see
Sections \ref{sec:collectively} and \ref{sec:Bier-compl})}.

\medskip
If $\mathcal{K} = \langle K_i\rangle_{i=1}^r = \langle K_1,\dots,
K_r\rangle$ is a collective $r$-unavoidable collection of
subcomplexes of $2^{[n]}$ (Definition~\ref{dfn:pigeonhole}), then
(by Problem~\ref{prob:main_problem}~(3)) the associated deleted
join $\mathcal{K}^\ast_\Delta = K_1\ast_\Delta\dots\ast_\Delta
K_r$ is expected to be $(n-r-1)$-connected. This is indeed the
case, as shown in the first part of Theorem~\ref{ThmMain}. In
particular we recover the result that $K\ast_\Delta K^\circ$ is an
$(n-2)$-dimensional homotopy sphere, whenever $K\neq 2^{[n]}$ is
{\em superdual} in the sense that $K\supseteq K^\circ$.

\medskip
In the special case when $\mathcal{K} = \langle
K_i\rangle_{i=1}^r$ is an Alexander $r$-tuple
(Definition~\ref{DfnAlexander}), we have a stronger result (see
the second half of Theorem~\ref{ThmMain}), that the associated
 Bier complex $\mathcal{K}^\ast_\Delta$ is a wedge of
spheres of the same dimension $n-r$ (Property ($3^\sharp$) in
Problem~\ref{prob:main_problem}). We describe an algorithm how the
number of these spheres can be explicitly calculated
(Corollary~\ref{cor:korolar}) and illustrate the calculation in
the case of `optimal chessboard complexes'
(Section~\ref{SecNumb}).

\medskip
A classification theorem for Alexander $r$-tuples
(Theorem~\ref{thm:class}) is proved in
Section~\ref{sec:classification}. It turns out, somewhat
unexpectedly and as a pleasant surprise, that the `optimal
chessboard complexes' (introduced in
Section~\ref{sec:multi-chess}) are the central examples of Bier
complexes (Section~\ref{sec:Bier-compl}) for $r\geq 3$.

%\medskip
%The isolation of Alexander $r$-tuples as potentially interesting
%objects of study lead to the natural question of finding
%interesting examples of these objects. As a pleasant surprise the
%`optimal chessboard complexes' (Section~\ref{sec:multi-chess})
%turned out to be  Bier complexes (Section~\ref{sec:collectively}).

\medskip
Among the corollaries of our results are exact connectivity bounds
for some classes of generalized chessboard complexes (including
the main case of Theorem~3.2 from \cite{jvz-1}). These results are
highly relevant for applications to the results of Tverberg-Van
Kampen-Flores type. As illustrated by the results in
Section~\ref{sec:collectively}, our alternative methods provide
some new insight complementing both the `constraint method' of
\cite{bfz} and the `index theory' approach \cite{Mat, jvz-3}.

\medskip The rest of the paper is organized as follows.
Section~\ref{secBasic} is an overview of basic notions and facts,
including a brief exposition of the discrete Morse theory \cite{Forman1,
Forman2} (which is our central tool in this paper). We
develop a version of this method which appears to be particularly
well adapted for the analysis of Bier spheres
(Section~\ref{sec:two-perfect}). We show in Section
\ref{SecUnavoid} how the method can be extended to the case of
deleted joins of collective $r$-unavoidable complexes and general
Bier complexes (introduced in Section \ref{sec:Bier-compl}). The
highlights include the construction of a perfect discrete Morse
function in the case of `optimal chessboard complexes'
(Section~\ref{SecNumb}) and their relatives `long chessboard
complexes' (Section~\ref{secLong}).

\subsection*{Acknowledgements}It is our pleasure to acknowledge the support and hospitality
of the {\em Mathematisches Forschungsinstitut Oberwolfach}, where
in the spring of 2016  this paper was initiated as a `research in
pairs' project. The construction from the proof of Theorem
\ref{ThmMain} and the constructions of Section
\ref{sec:two-perfect} are supported by the Russian Science
Foundation under grant 16-11-10039. R. \v{Z}ivaljevi\'{c}
acknowledges the support of the Ministry of Education, Science and
Technological Development of Republic of Serbia, Grant 174034.

\section{An overview of basic definitions and facts}\label{secBasic}

In this section we collect some standard definitions and facts, as
a reminder for the reader. This is also an opportunity to
introduce some less standard notation and concepts, used in the
rest of the paper. For other standard facts and definitions the
reader is referred to \cite{Mat}.

%\textbf{(HERE STARTS THE PAPER:)}

\subsection{Simplicial complexes}

A simplicial complex on a set $V$ of vertices is a subset
$K\subset 2^V$ such that (1) $\emptyset\in K$ and (2) if $A\subset
B\in K$ then $A\in K$. By definition it is possible that
$\{v\}\notin K$ for some $v\in V$, however $K\neq\emptyset$ (since
$\emptyset\in K$ by property (1)).

\medskip
The complex $2^V$ is often referred to as the simplex spanned by
$V$, and denoted by $\Delta(V)$. We use, side by side, topological
and combinatorial language (and notation). For example,
$${[n] \choose \leqslant k},$$
is the $(k-1)$-skeleton of the $(n-1)$-dimensional simplex
$\Delta([n])$.

\medskip
The {\em deleted join} \cite[Section~6]{Mat} of a family
$\mathcal{K} = \langle K_i\rangle_{i=1}^r = \langle K_1,\dots,
K_r\rangle$ of subcomplexes of $2^{[n]}$ is the complex
$\mathcal{K}^\ast_\Delta = K_1\ast_\Delta\dots \ast_\Delta  K_r
\subset (2^{[n]})^{\ast r}$ where $A = A_1\uplus\dots\uplus A_r\in
\mathcal{K}^\ast_\Delta$ if and only if $A_j$ are pairwise
disjoint and $A_i\in K_i$ for each $i=1,\dots, r$.

\subsection{Multiple chessboard complexes}\label{sec:multi-chess}

A `chessboard complex', in a very broad sense, is any subcomplex
$K\subset 2^{([n]\times [r])}$ of the simplex $\Delta([n]\times
[r])$ spanned by elementary squares of an $(n\times
r)$-chessboard. Following \cite[Section~2.1]{jvz-1}, the multiple
chessboard complex
$$\Delta_{n,r}^{m_1,\dots, m_r; \mathbf{1}} =
\Delta_{n,r}^{m_1,\dots, m_r; 1,\dots, 1}$$ is described by the
condition that $S\in \Delta_{n,r}^{m_1,\dots, m_r; \mathbf{1}}$ if
and only if the cardinality of the set $S\cap ([n]\times \{i\})$
is at most $m_i$ for each $i=1,\dots, r$, and the cardinality of
the set $S\cap (\{j\}\times [r])$ is at most $1$ for each
$j=1,\dots, n$.

A moment's reflections reveals that there is a relation,
$$ \Delta_{n,r}^{m_1,\dots, m_r; \mathbf{1}} \cong
{[n]\choose \leqslant m_1}\ast_\Delta\dots\ast_\Delta {[n]\choose
\leqslant m_r},$$ which says that the multiple chessboard complex
can be always expressed as the deleted join of skeletons of the
simplex $\Delta([n])\cong \Delta^{n-1}$.

\medskip
One of the central results of \cite[Theorem~3.2]{jvz-1} says that
$\Delta_{n,r}^{m_1,\dots, m_r; \mathbf{1}}$ is $(\nu-2)$-connected
where $\nu = m_1+\dots +m_r$, provided $n \geq m_1+\dots +m_r +
r-1$. For this reason the chessboard complex
$\Delta_{n,r}^{m_1,\dots, m_r; \mathbf{1}}$ is often called {\em
optimal}, if $n = m_1+\dots +m_r + r-1$.
Similarly we say that a multiple chessboard complex is {\em long}
if $n > m_1+\dots +m_r + r-1$.

\subsection{ Alexander duality and Bier spheres}\label{section_pairings}

The Alexander dual \cite[Section~5.6]{Mat} of $K\subset 2^V$ is
the set $K^\circ$ of all complements of non-simplices in $K$,
$$K^\circ = \{F\subset V \mid V\setminus F\notin K\}.$$
In order to rule out the possibility $K^\circ =
\emptyset$, we tacitly assume throughout the paper that $K\neq
2^V$, whenever we are dealing with Alexander pairs $(K, K^\circ)$
of complexes.

For a given simplicial complex $K \subset 2^{[n]}$, the associated
Bier sphere $Bier(K) = K \ast_{\Delta} K^{\circ}$ is described as
the deleted join of $K$ with its Alexander dual  $K^{\circ}$. The
simplices of the deleted join $K\ast_{\Delta}K^{\circ}$ are by
definition disjoint unions $A_1\uplus A_2\subset [n]\uplus [n]
\cong [n]\times [2]$, where $A_1\in K, A_2\in K^\circ$ and
$A_1\cup A_2\neq\emptyset$. They can be also described as ordered
partitions of the set $[n]$ into three parts $(A_1,A_2;B)$ (where
$B:=[n]\setminus (A_1\cup A_2)$).

\medskip

Note that a partition $(A_1,A_2;B)$ corresponds to a simplex in
the deleted join $K\ast_\Delta K^\circ$ if and only if:
\begin{enumerate}
 \item $A_1\in K$,
 \item $A_2\in K^{\circ}$ (or equivalently $[n]\setminus A_2
 \notin K$);
 \item $\emptyset\neq B\neq [n]$ (equivalently $\emptyset \neq A_1\cup A_2\neq
 [n]$).
\end{enumerate}

The incidence relation of the simplices is described by the rule:

$(A_1,A_2;B) \subseteq (A'_1,A_2';B')$ iff $A_1\subseteq A'_1$, and
$A_2\subseteq A_2'$.

\subsection{Discrete Morse theory}\label{section_prelim}

Robin Forman's discrete Morse theory \cite{Forman1, Forman2} is, as a tool,  as powerful as the smooth
Morse theory. It has been used in computations of the homology, the cup-product, Novikov homology, and other
topological  and combinatorial computations and applications.
Major advantage of discrete Morse theory (compared to smooth Morse theory) is its applicability
to a considerably larger class of objects which include simplicial and cellular complexes (and not only smooth manifolds).

\smallskip
In our paper we make use of a relatively small and quite reduced piece of the general theory. For our purposes
it suffices   to think of a  `Morse function' as a special kind of
matching on the set of simplices. Here is a brief overview of some of the central definitions and results of discrete
Morse theory.

\medskip

Let $K$ be a simplicial complex. Its $p$-dimensional simplices
($p$-simplices for short) are denoted by $\alpha^p, \ \beta^p,
\sigma^p$, etc. A \textit{discrete vector field} $D$ is a set of
pairs $\big(\alpha^p,\beta^{p+1}\big)$ (called a matching) such
that:
\begin{enumerate}
    \item  each simplex of the complex participates in at most one
    pair, and
    \item  in each pair, the simplex $\alpha^p$ is a facet of $\beta^{p+1}$.
\end{enumerate}
\noindent The pair $(\alpha^p, \beta^{p+1})$ can be informally
thought of as a vector in the vector field $D$. For this reason it
is occasionally denoted by $\alpha^p \rightarrow \beta^{p+1}$ (and
in this case $\beta^{p+1}$ is referred to as  {\em the end}\/ of
the arrow $\alpha^p \rightarrow \beta^{p+1}$).

\smallskip
Given a discrete vector field $D$, a \textit{gradient path} in $D$
is a sequence of simplices

$$\alpha_0^p, \ \beta_0^{p+1},\ \alpha_1^p,\ \beta_1^{p+1}, \ \alpha_2^p,\ \beta_2^{p+1} ,..., \alpha_m^p,\ \beta_m^{p+1},\ \alpha_{m+1}^p,$$
which satisfies the following conditions:
\begin{enumerate}
\item $p\geq 0$, that is, the empty set $\emptyset\in K$ is never
matched,
    \item  $\big(\alpha_i^p,\ \beta_i^{p+1}\big)$ is a pair in $D$ for each $i$,
    \item for each $i = 0,\dots, m$ the simplex $\alpha_{i+1}^p$
    is a facet of $\beta_i^{p+1}$.
    \item $\alpha_i\neq \alpha_{i+1}$.
\end{enumerate}

A path is \textit{closed} if $\alpha_{m+1}^p=\alpha_{0}^p$. A
\textit{discrete Morse function } (DMF for short) is a discrete
vector field without closed paths.

\medskip
Assuming that a discrete Morse function is fixed, the {\em
critical simplices} are those simplices of the complex that are
not matched. The Morse inequality \cite{Forman2} states that critical
simplices cannot be completely avoided.

A discrete Morse function is a  \textit{perfect Morse function}
whenever the number of critical $k$-simplices  equals the $k$-th
Betty number of the complex.  It is equivalent to the condition
that the number of all critical simplices equals the sum of Betty
numbers.

\medskip
Perhaps the main idea of discrete Morse theory, as summarized in
the following theorem of R.~Forman, is to contract all matched
pairs of simplices and to reduce the simplicial complex $K$ to a
cell complex (where critical simplices correspond to the cells).

\begin{thm}\cite{Forman1, Forman2} \label{ThmWedge} Assume that a
discrete Morse function on a simplicial complex $K$ has a single
zero-dimensional critical simplex $\sigma^0$ and that all other
critical simplices have the same dimension $N>1$. Then $K$ is
homotopy equivalent to a wedge of $N$-dimensional spheres.

 More
generally, if all critical simplices, aside from $\sigma^0$, have
dimension $\geq N$, then the complex $K$ is $(N-1)$-connected.  \qed
\end{thm}

\subsection{The `constraint method' and `unavoidable complexes'}
\label{sec:constraint}

The Gro\-mov-Blagojevi\'{c}-Frick-Ziegler reduction, or the {\em
constraint method}, is an elegant and powerful method for proving
results of Tverberg-Van Kampen-Flores type. It relies on the
concept of `unavoidable' or more precisely $r$-unavoidable
complex, where $r$ is a positive integer. The property of being
`unavoidable' is one of the central themes of our paper. For this
reason we briefly review the `constraint method' where this
concept originally appeared.

\begin{equation}
\begin{CD}
K @>f>> \mathbb{R}^d\\
@VeVV @ViVV\\
\Delta^N @>F>> \mathbb{R}^{d+1}
\end{CD}
\end{equation}

\noindent Suppose that the continuous Tverberg theorem holds for
the triple $(\Delta^N, r, \mathbb{R}^{d+1})$ in the sense that for
each continuous map $F: \Delta^N\rightarrow \mathbb{R}^{d+1}$
there exists a collection of $r$ vertex disjoint simplices
$\Delta_1,\ldots, \Delta_r$ of $\Delta^N$ such that
$F(\Delta_1)\cap\ldots\cap F(\Delta_r)\neq\emptyset$. For example
the {\em Topological Tverberg theorem} \cite[Section~6]{Mat}
(proved by B\'{a}r\'{a}ny, Shlosman, and Sz\"{u}s for primes, and
\"{O}zaydin for prime powers) says that this is the case if $r =
p^k$ is a prime power and $N = (r-1)(d+2)$. Suppose that $K\subset
\Delta^N$ is a simplicial complex which is {\em $r$-unavoidable}
in the sense that if $A_1\uplus\ldots\uplus A_r = [N+1]$ is a
partition of the set $[N+1]$ (of vertices of $\Delta$), then at
least one of the simplices $A_i$ of $\Delta^N$ is in $K$. Then for
each continuous map $f : K\rightarrow \mathbb{R}^d$ there exists
vertex disjoint simplices $\sigma_1,\ldots, \sigma_r\in K$ such
that $f(\sigma_1)\cap\ldots\cap f(\sigma_r)\neq\emptyset$.

\medskip
Indeed, let $F'$ be an extension ($F'\circ e = f$) of the map $f$
to $\Delta^N$. Suppose that $\rho : \Delta^N\rightarrow
\mathbb{R}$ is the function $\rho(x) := {\rm dist}(x, K)$,
measuring the distance of the point $x\in \Delta^N$ from $K$.
Define $F = (F', \rho)  : \Delta^N\rightarrow \mathbb{R}^{d+1}$
and assume that $\Delta_1,\dots, \Delta_r$ is the associated
family of vertex disjoint simplices of $\Delta^N$, such that
$F(\Delta_1)\cap\ldots\cap F(\Delta_r)\neq\emptyset$. More
explicitly suppose that $x_i\in\Delta_i$ such that $F(x_i)=F(x_j)$
for each $i,j = 1,\ldots, r$. Since $K$ is $r$-unavoidable,
$\Delta_i\in K$ for some $i$. As a consequence $\rho(x_i)=0$, and
in turn $\rho(x_j)=0$ for each $j=1,\ldots, r$. If $\Delta_i'$ is
the minimal simplex of $\Delta^N$ containing $x_i$ then
$\Delta_i'\in K$ for each $i=1,\ldots, r$ and
$f(\Delta_1')\cap\ldots\cap f(\Delta_r')\neq\emptyset$.

\medskip
The reader is referred to \cite{bfz} for a more complete
exposition and numerous examples of applications of the
`constraint method',  see also \cite[Section~2.9(c)]{Grom-2} for
the historically first appearance of the idea.

\section{Collectively unavoidable $r$-tuples\\ and Alexander $r$-tuples  of complexes}
\label{sec:collectively}

In this section we introduce the central objects of our paper. Our
tacit assumption is that all complexes $K$ are {\em proper
subcomplexes} of $2^{[n]}$ in the sense that $K\subsetneq
2^{[n]}$.

\begin{dfn}\label{dfn:pigeonhole}
An ordered $r$-tuple $\mathcal{K} = \langle K_1,...,K_r\rangle$ of
subcomplexes of $2^{[n]}$ is {\em collective $r$-unavo\-id\-able}
(we also say that $\mathcal{K}$ is a {\em pigeonhole $r$-tuple} on
$[n]$), if for each ordered collection $(A_{1},...,A_{r})$ of
disjoint sets in $[n]$ there exists $i$ such that $A_i\in K_i$.
The class of collective $r$-unavoidable complexes is denoted by
$CU_{r,n}$, or by $CU_r$ if $n$ is fixed or clear from the
context.
\end{dfn}

On closer inspection,  the definition can be usefully rephrased as
follows.  For the ordered $r$-tuple $\mathcal{K} = \langle
K_{1},\dots, K_{r}\rangle$ and for an ordered disjoint collection
$(A_{1}, \dots, A_{r})$ of subsets of $[n]$, we construct a
bipartite graph $\Gamma \subset K_{r,r}$, where by definition
there is an edge $(i,j)\in\Gamma$ if and only if $A_{i}\notin
K_{j}$. Then the collective $r$-unavoidability of $\mathcal{K}$ is
equivalent to the condition that the graph $\Gamma$ does not
contain a complete matching (does not satisfy the {\em marriage
condition} of the classical Hall's theorem). We therefore conclude
that the pigeonhole property does not depend on the ordering of
simplicial complexes.

\begin{rem} The bipartite graph $\Gamma = \{(i,j)\in [r]^2 \mid A_i\notin
K_j\}$ interpretation naturally leads to an extension of
Definition~\ref{dfn:pigeonhole} to the case of collections
$\mathcal{K} = \langle K_1,...,K_s\rangle$ where $s$ is not
necessarily equal to $r$. Note however that the symmetric case
$s=r$ is somewhat exceptional. For example the classical
`Hilfssatz' of Frobenius \cite{Sch} implies that $\mathcal{K} =
\langle K_1,...,K_r\rangle$ is collective $r$-unavoidable if and
only if for each ordered collection $(A_{1},...,A_{r})$ of
disjoint sets in $[n]$ there exists a pair $(S,T)$ of subsets of
$[r]$, such that $\vert S\vert +\vert T\vert = r+1$, and $A_i\in
K_j$ for each $i\in S$ and $j\in T$.
\end{rem}

It is easy to characterize all pigeonhole $2$-tuples: $(K_{1},
K_{2})$ is collective unavoidable if and only if
$K_{1}^{\circ}\subset K_{2}$ (or equivalently
$K_{2}^{\circ}\subset K_{1}$).

For an $r$-tuple of complexes $\langle K_1,...,K_r\rangle$ we
shall use a natural partial ordering on the set of all set of
pairwise disjoint $r$-tuples $(A_1,...,A_r)$ with $A_i\in K_i$:
say that $(A_1,...,A_r)\leqslant(A'_1,...,A'_r) $   whenever
$\forall i: A_i\subseteq A'_i$.

We also put a partial ordering on the set of all $r$-tuples of
complexes by the same rule. So we automatically have the notion
of\textit{ minimal unavoidable $r$-tuple of complexes} $\langle
K_1,...,K_r\rangle$.

\begin{lemma}\label{LemmaUnavoid} Suppose that the $r$-tuple
$\mathcal{K} = \langle K_1,...,K_r\rangle$ is \textit{collective
$r$-unavoidable}. Then for each maximal disjoint collection
$(A_{1},...,A_{r})$  with $A_{i}\in K_i$, the set $[n]\setminus
\bigcup_{i=1}^{r} A_i$ contains at most $r-1$ elements.
\end{lemma}

\textit{Proof.} Suppose that $\langle K_1,...,K_r\rangle$ is
collective $r$-unavoidable. Let $(A_{1},...,A_{r})$ be a maximal
disjoint collection satisfying the condition $A_{i}\in K_i$ for
each $i=1,\dots ,r$. Suppose (for contradiction) that
$\{a_1,\ldots, a_r\}\subset [n]\setminus \bigcup_{i=1}^{r} A_i$.
Then $A_i'=A_i\cup \{a_i\}\notin K_i$ (by the maximality of the
collection $(A_{1},...,A_{r})$) and the collection
$(A_{1}',...,A_{r}')$ clearly violates the collective
$r$-unavoidability condition for $\langle K_1,...,K_r\rangle$.
\qed

\begin{dfn}\label{DfnAlexander}
An $r$-tuple of complexes $\mathcal{K} = \langle
K_1,...,K_r\rangle$  on one and the same set of vertices $[n]$ is
an \textit{Alexander $r$-tuple} if,
\begin{enumerate}
\item it is collective $r$-unavoidable, and \item for each
$r$-tuple of sets $A_1,...,A_r$ with $A_i \in K_i$ the set
$[n]\setminus \bigcup_{i=1}^{r}A_i$ has at least $r-1$ elements.
\end{enumerate}
The class of Alexander $r$-tuples of subcomplexes of $2^{[n]}$ is
denoted by $\mathcal{A}_{r}$ (or by $\mathcal{A}_{r,n}$ if the set
$[n]$ of vertices should be emphasized).
\end{dfn}

\begin{prop}\label{prop:Pure} Given an Alexander $r$-tuple on $[n]$, for each maximal
$r$-tuple of disjoint sets $(A_1,...,A_r)$ with $A_i \in K_i$ the
set $[n]\setminus \bigcup_{i=1}^{r}A_i$ has exactly  $r-1$
elements.
\end{prop}
\textit{Proof.} This follows from Lemma \ref{LemmaUnavoid} and the
property (2) from the definition of the Alexander $r$-tuple
(Definition~\ref{DfnAlexander}).\qed

\begin{prop}\label{prop:Alex-min}
An Alexander $r$-tuple of complexes is always a minimal pigeonhole
$r$-tuple of complexes.
\end{prop}

\textit{Proof.} Assume $\langle K_1,\dots,K_r \rangle$ is an
Alexander $r$-tuple which is not a minimal collective
$r$-unavoidable collection of complexes. This means that (possibly
after a re-enumeration) the collection $\langle K_1\setminus
\{A_1\}, K_2,\dots,K_r\rangle$ is also collective $r$-unavoidable
for some maximal simplex $A_1\in K_1$. As a consequence the
restrictions $\langle K_2|_{[n]\setminus
A_1},...,K_r|_{[n]\setminus A_1}\rangle$ form a collective
$(r-1)$-unavoidable family of complexes. Lemma \ref{LemmaUnavoid}
implies that for any  maximal disjoint collection $(A_2,...,A_r)$
such that $A_j\in K_j|_{[n]\setminus A}$ for each $j=2,\dots, r$,
the set $[n]\setminus \bigcup_{i=1}^r A_{i}$ contains strictly
less than $r-1$ elements. Then $(A_1,\dots,A_r)$ is a maximal
family satisfying $A_i\in K_i$ for each $i=1,\dots, r$, which is
in contradiction with the condition (2) from
Definition~\ref{DfnAlexander}. \qed

\bigskip

The converse of Proposition~\ref{prop:Alex-min} is in general not
true.

\begin{Ex}
 $$K_1=K_2=K_3={[10]\choose \leqslant 2} \sqcup {[9]\choose \leqslant 3} $$
is a minimal collective unavoidable  $3$-tuple which is not an
Alexander $3$-tuple.
\end{Ex}

\begin{Ex}\label{ex:Alex-para} A $2$-tuple of complexes is
an Alexander $2$-tuple iff it is a pair of mutually dual complexes
$(K,K^{\circ})$.
\end{Ex}

 The following examples describes the Alexander
$r$-tuples $\mathcal{K} =  \langle K_1,\dots, K_r\rangle$ where
each of the complexes $K_i$ is a skeleton of the simplex
$2^{[n]}$.

\begin{Ex}\label{ex:skeletons}
The collection of subcomplexes of $2^{[n]}$,
$$\left({[n]\choose \leqslant m_1},\dots,{[n]\choose \leqslant
m_r}\right)$$ is always an Alexander $r$-tuple, provided
$n=\sum_{i=1}^r m_i+r-1$.
\end{Ex}

\begin{Ex}\label{ex:Alex-ex}
Define a simplicial complex $K\subset 2^{[6]}$ as the cone with
apex $1$ over the five-element set $\{2,3,4,5,6\}$. The complex
$K$ is essentially a graph with five edges
$\{1,2\},\{1,3\},\{1,4\},\{1,5\},\{1,6\}$. It is not difficult to
see that $\langle K,K,K\rangle$ is indeed an Alexander $3$-tuple.

\smallskip This example is the simplest case of a more general
construction. For given integers $m_1, m_2,\ldots, m_r$, let $n =
m_1+m_2+\cdots+m_r+r-1$. Choose a simplex $\Delta(C)$ spanned by
$C\neq\emptyset$ (where $[n]\cap C = \emptyset$) and define the
complexes,
$$K_i={[n]\choose \leqslant m_i}*\Delta(C)\textrm{ , for all }i=1,2,\ldots,r.$$
It can be easily seen that $\langle K_1,K_2,\ldots,K_r\rangle$ is
an Alexander $r$-tuple.
\end{Ex}

\subsection{Operations generating collective $r$-unavoidable
complexes}

As de\-mon\-strated by the classification theorem
(Theorem~\ref{thm:class}), Alex\-an\-der $r$-tuples are scarce,
and a very special class of simplicial complexes. The situation
with the collective $r$-unavoidable complexes is quite the
opposite, as illustrated by the following construction.

\smallskip Let $r\geq 2$ and let $K_i\subset 2^{[n]}$ be a
collection of not necessarily distinct simplicial complexes.
Assume that the $(r-1)$-tuple $\mathcal{K} = \langle
K_1,K_2,\ldots,K_{r-1}\rangle$ is NOT collective
$(r-1)$-unavoidable on $[n]$.

\smallskip
Define $R(\mathcal{K}) = R_r(\mathcal{K})=
R_r(K_1,K_2,\ldots,K_{r-1})$ as the  subcomplex of $2^{[n]}$ where
$F\in R(\mathcal{K})$ if and only if there exists an ordered
partition $F_1\uplus\dots\uplus F_{r-1} = F^c$ of the complement
of $F$ such that $F_i\notin K_i$ for each $i=1,\dots, r-1$.

\medskip Observe that $R(\mathcal{K})$ is generated by the sets $(F_1\cup\dots\cup
F_{r-1})^c$ where $F_j$ are pairwise disjoint and $F_j$ is a
minimal non-face in $K_j$ for each $j=1,\dots, r-1$.

\medskip Note that $R_r(\mathcal{K})$ can be described as the unique minimal
simplicial complex $Z$ such that $\langle
K_1,\ldots,K_{r-1},Z\rangle$ is a collective unavoidable $r$-tuple
on $[n]$. Observe that $\emptyset\in R_r(\mathcal{K})$ follows
from the assumption that $\mathcal{K}$ is not $(r-1)$-unavoidable.

\begin{dfn}
The complex $R_r(\mathcal{K})$ is referred to as the {\em residual
complex} of the $(r-1)$-tuple $\mathcal{K} = \langle
K_1,K_2,\ldots,K_{r-1}\rangle$.
\end{dfn}
Observe that in the case $r=2$ the residual complex of $K\subset
2^{[n]}$ is precisely the Alexander dual, $R(K)=K^{\circ}$. More
generally, for a complex $K\subseteq 2^{[n]}$ we define the
associated $r^{\rm th}$ residual complex
$R_r(K)=R(K_1,\ldots,K_{r-1})$ where $K_1=\dots=K_{r-1}=K$. Note
that $K$ is a minimal $r$-unavoidable complex if and only if
$R_r(K)=K$.

\begin{prob}
Find interesting examples of ordered collections of complexes
$\mathcal{K} = \langle K_1,K_2,\ldots,K_{r-1}\rangle$ such that
$\langle K_1,K_2,\ldots,K_{r-1}, R(\mathcal{K})\rangle$ satisfies
the condition $(3^\sharp)$ (in Problem~\ref{prob:main_problem}).
\end{prob}

\section{Bier complexes}\label{sec:Bier-compl}

For each Alexander $2$-tuple $\langle K_1, K_2\rangle = \langle
K_1, K_1^\circ\rangle = \langle K_2^\circ, K_2\rangle$, the
associated deleted join $K_1\ast_\Delta K_2$ is the standard Bier
sphere $Bier(K_1)\cong Bier(K_2)$ (Example~\ref{ex:Alex-para}).
This observation motivates the following definition.

\begin{dfn}\label{dfn:r-Bier} Suppose that $\mathcal{K} = \langle K_1,\dots, K_r\rangle$
is an Alexander $r$-tuple of complexes $K_i\subset 2^{[n]}$. Then
the associated {\em  Bier complex} is defined as the
deleted join,
\[ Bier(\mathcal{K}) := \mathcal{K}^{\ast
r}_\Delta = K_{1}*_{\Delta}K_{2}*_{\Delta}...*_{\Delta}K_{r}.
\]
\end{dfn}

It is well known that the `join' and the `deleted join' operations
commute (see for example Lemma~6.4.3. in \cite{Mat}). The
following lemma is a natural generalization.

\begin{lemma}\label{lemma:commute}
Let $\mathcal{K} = \langle K_1,\dots, K_r\rangle$ and $\mathcal{L}
= \langle L_1,\dots, L_r\rangle$ be two collections of simplicial
complexes where $K_i\subset 2^{[m]}$ and $L_i\subset 2^{[n]}$ for
each $i=1,\dots, r$. Then,
\begin{equation}\label{eqn:commute}
(\mathcal{K}\ast \mathcal{L})^{\ast r}_\Delta \cong
\mathcal{K}^{\ast r}_\Delta \ast \mathcal{L}^{\ast r}_\Delta
\end{equation}
where by definition $\mathcal{K}\ast \mathcal{L} :=  \langle
K_1\ast L_1,\dots, K_r\ast L_r\rangle$.
\end{lemma}

\textit{Proof.} If $A = A_1\uplus\dots\uplus A_r\in
\mathcal{K}^{\ast r}_\Delta$ and $B= B_1\uplus\dots\uplus B_r\in
\mathcal{L}^{\ast r}_\Delta$ then $A\ast B\in \mathcal{K}^{\ast
r}_\Delta \ast \mathcal{L}^{\ast r}_\Delta$ corresponds to the
simplex $C_1\uplus\dots\uplus C_r\in (\mathcal{K}\ast
\mathcal{L})^{\ast r}_\Delta$ where $C_i := A_i\uplus B_i$ for
each $i=1,\dots, r$. \hfill $\square$

\bigskip The following theorem is one of the main results of our
paper. It says that the classes $CU_r$ and $\mathcal{A}_r$
respectively, satisfy the central properties $(3)$ and
$(3^\sharp)$, listed in Problem~\ref{prob:main_problem}.

\begin{thm}\label{ThmMain} Let $\mathcal{K} = \langle K_1,\dots,
K_r\rangle$ be a collection of subcomplexes of $2^{[n]}$.

\begin{enumerate}
 \item The deleted join $\mathcal{K}^{\ast r}_\Delta =
K_{1}*_{\Delta}K_{2}*_{\Delta}...*_{\Delta}K_{r}$ of a collective
$r$-unavoidable collection $\mathcal{K}$ of complexes is always
$(n-r-1)$-connected.
 \item The  Bier complex
$Bier(\mathcal{K}) =
K_{1}*_{\Delta}K_{2}*_{\Delta}...*_{\Delta}K_{r}$, associated to
an Alexander $r$-tuple $\mathcal{K}$, is a pure complex of
dimension $n-r$, homotopy equivalent to a wedge of
$(n-r)$-dimensional spheres.
 \end{enumerate}
 \end{thm}

The following corollary of the proof of Theorem~\ref{ThmMain}
emphasizes the computational efficiency of the approach based on
the discrete Morse function described in Section~\ref{SecUnavoid}.

\begin{cor}\label{cor:korolar}
For an Alexander $r$-tuple $\mathcal{K}$ the number of spheres in
the wedge $Bier(\mathcal{K})$ can be efficiently calculated as the
number of critical simplices of the discrete Morse function $D$
constructed in Section \ref{SecUnavoid}.
\end{cor}

The efficiency of the method is illustrated in
Section~\ref{SecNumb} by the calculation of the number of spheres
in the important particular case of the optimal multiple
chessboard complex (Section~\ref{sec:multi-chess}).

Recall that the number of spheres in a wedge decomposition can be
in principle calculated as the reduced Euler characteristic of the
complex. This calculation is typically very slow and inefficient,
as it is based on an `inclusion-exclusion' type formula which
involves enumeration of all simplices in $Bier(\mathcal{K})$.

\medskip
One of important motivations for introducing (collective)
$r$-unavoidable complexes are applications to problems of
Tverberg-Van Kampen-Flores type. By emphasizing the role of
Theorem~\ref{ThmMain}, the following corollaries provide some
initial evidence illustrating this interesting and important
connection.

\begin{cor} \label{cor:index-inequality}
{\rm (\cite[Theorem~5.5.5]{Mat}, \cite[Theorem~3.6]{jvz-3})}
Suppose that $K$ is an $r$-unavoidable complex with vertices in
$[n]$. Suppose that $r=p^k$ is a prime power and let $G =
(\mathbb{Z}_p)^k$ be an elementary abelian $p$-group acting freely
on the set $[r]$. Let $K^{\ast r}_\Delta$ be the $r$-fold deleted
join of $K$. Then,
\begin{equation}\label{eqn:izvor}
{\rm Ind}_G (K^{\ast r}_\Delta) \geq n - r,
\end{equation}
where ${\rm Ind}_G$ is the associated equivariant index function
\cite{Mat,jvz-3}.
\end{cor}

\textit{Proof.} If $K\subset 2^{[n]}$ is $r$-unavoidable then the
collection $\mathcal{K} = \langle K_1,\dots, K_r \rangle$, where
$K_1=\dots = K_r = K$, is a collective $r$-unavoidable collection
of complexes. By Theorem~\ref{ThmMain} the deleted join $K^{\ast
r}_\Delta$ is $(n-r-1)$-connected. The inequality
(\ref{eqn:izvor}) follows from this observation and the basic
properties of the index function ${\rm Ind}_G$, see for
example Proposition~3.3 (inequality (5)) in \cite{jvz-3}. \qed

\medskip
The following result is a simplest example illustrating the role
of $r$-unavoid\-able complexes in Tverberg type problems. For a
more general theorem of this type the reader is referred to
\cite[Theorem~4.4]{bfz}, see also \cite[Theorem~4.6]{jvz-3} for a
related result.

\begin{cor} {\rm (\cite{bfz})}
Suppose that $K\subset 2^{[n]}$ is an $r$-unavoidable complex.
Assume that $r=p^k$ is a prime power and let $d$ be the integer
satisfying the inequality $(r-1)(d+2)+1\leq n$. Then $K$ is {\em
globally $r$-non-embeddable in $\mathbb{R}^d$} in the sense that
for each continuous map $f : K\rightarrow \mathbb{R}^d$ there
exist $r$ vertex-disjoint simplices $\Delta_1,\dots, \Delta_r$ of
$K$ such that,
\[
f(\Delta_1)\cap\dots\cap f(\Delta_r)\neq\emptyset.
\]
\end{cor}

\textit{Proof.}  The most elegant proof of this result is by the
`constraint method' \cite{bfz} (see Section~\ref{sec:constraint}
for an outline). The `index theory proof', in the spirit of
\cite[Section~6]{Mat} and \cite{jvz-3}, is based on
Corollary~\ref{cor:index-inequality}.          \qed

\medskip\noindent
\begin{rem}{\rm Let us observe that the
`Gromov-Blagojevi\'{c}-Frick-Ziegler reduction' (the `constrain
method') reduces a Van Kampen-Flores (or Tverberg) type question,
to another result of that type. More explicitly (and more
generally) the method says that the question if there exists a map
$f : K\rightarrow \mathbb{R}^d$ without (global) $r$-fold points
(Tverberg $r$-tuples) can be reduced to a similar problem for an
appropriate map $F : \Sigma \rightarrow \mathbb{R}^D$. Here
$K\subset\Sigma$ is a complex which is {\em relatively
$r$-unavoidable subcomplex of $\Sigma$} in the sense of
\cite[Definition~2.5]{jvz-3}.

This reasoning illustrates why the `index theory methods' (which
rely on results of Dold and Volovikov, see
\cite[Section~6.2.6]{Mat}) retain their importance. This also
explains why the results like Theorem~\ref{ThmMain} may be
interesting since both the Dold's and the Volovikov's theorem are
based on the homotopical (respectively homological) connectivity
of the associated configuration space (deleted join).

For illustration, Theorem~2.1 from \cite{jvz-2}, that needs such a
connectivity result for its proof, is possibly a good candidate
for a Tverberg-Van Kampen-Flores type result that cannot be
obtained directly by the `constraint method'. }
\end{rem}

\subsection{Bier complexes and discrete Morse theory}

The proof of Theorem~\ref{ThmMain} (Section~\ref{SecUnavoid}) and
the proofs of other connectivity results in this paper rely on
Discrete Morse theory (Theorem~\ref{ThmWedge}). All our discrete
Morse functions (DMF) are defined on deleted joins
$\mathcal{K}^{\ast r}_\Delta = K_1\ast_\Delta\dots\ast_\Delta K_r$
of complexes $K_i\subset 2^{[n]}$ and they all have some common
features.

\medskip A simplex $\beta \in \mathcal{K}^{\ast r}_\Delta$ is
usually recorded as a disjoint sum $\beta = A_1\uplus\dots\uplus
A_r$, see \cite[Sections~5 and 6]{Mat}. We find it convenient (for
bookkeeping purposes) to use an alternative `partition notation'
$\beta = (A_1,\dots, A_r; B)$ where $B = [n]\setminus
\cup_{i=1}^r~A_i$. To match a $p$-simplex $\alpha^p = (A_1',\dots,
A_r'; B')$ with a $(p+1)$-simplex $\beta^{p+1} = (A_1,\dots, A_r;
B)$ is the same as to {\em `migrate'} an element $i\in B'$ to one
of the sets $A_j'$. This is possible if $\alpha^p$ is a facet of
$\beta^{p+1}$ i.e.\ if $B' = B\uplus \{i\}$ for some $i\in B'$.

\medskip\noindent {\bf Caveat:} { In the paper}, we simplify the
notation by omitting the braces and by writing simply $B\cup i$
instead of $B\cup \{i\}$ (with the tacit assumption that $i\notin
B$). We also write $j < B$ ($j>B$) if $j<i$ for each $i\in B$
(respectively if $j>i$ for each $i\in B$).

\section{Classification theorem for Alexander $r$-tuples}
\label{sec:classification}

\begin{thm}\label{thm:class} If $\mathcal{K} = \langle K_1,\dots, K_r\rangle$
is an Alexander $r$-tuple then,
\begin{enumerate}
 \item[(1)]\hspace{0.5mm} $r=2$ and $(K_1, K_2) = (K, K^\circ)$ is an Alexander
 pair of dual complexes (Example~\ref{ex:Alex-para}), or
 \item[(2)] \hspace{0.5mm} $r\geq 3$ and $K_i ={[n] \choose \leq m_i
 }$ (Example~\ref{ex:skeletons}) where $n = m_1+\dots +m_r
 +r-1$, or
 \item[(3)] \hspace{0.5mm} $r\geq 3$ and $K_i ={[n] \choose \leq m_i  }\ast
 \Delta(C)$ (Example~\ref{ex:Alex-ex}) where $n = m_1+\dots +m_r
 +r-1$ and $\Delta(C) = 2^C$ is the simplex spanned by a non-empty set $C$ such that $C\cap
 [n]=\emptyset$.
\end{enumerate}
\end{thm}

\textit{Proof.} Suppose that $r\geq 3$. A minimal non-simplex of a
simplicial complex $K\subset 2^{[n]}$ is called a $K$-blocker.
Equivalently, $A\subset [n]$ is a $K$-blocker if $A\notin K$ and
$\partial(A)\subset K$.

\smallskip
Let $\mathcal{A} = (A_1,\dots, A_r)$ be a maximal disjoint
$r$-tuple of sets in $[n]$ such that $A_i\in K_i$ for each
$i=1,\dots, r$. Moreover, we assume that $A_r$ has the maximal
size possible in all such $r$-tuples.

\smallskip

Since $\mathcal{K}$ is an Alexander $r$-tuple the set
$[n]\setminus \bigcup_{i=1}^r~A_i = \{t_1,\dots, t_{r-1}\}$ has
exactly $(r-1)$ elements (Proposition~\ref{prop:Pure}).

\smallskip
Let $X_{\mathcal{A}} = X_1\uplus\dots\uplus X_{r-1}\uplus X_r$ be
the associated `blocker partition' where $X_i := A_i\cup\{t_i\}$
for each $i=1,\dots, r-1$ and $X_r := A_r =
[n]\setminus\bigcup_{i=1}^{r-1}~X_i$. The name is justified by the
fact that $X_i$ is  a $K_i$-blocker for each $i=1,\dots, r-1$.

Indeed, suppose that $X_\nu$ is not a $K_\nu$-blocker for some
$\nu = 1,\ldots, r-1$, which means that there exists $x\in A_\nu$
such that $X_\nu\setminus\{x\}\notin K_\nu$. The maximality of
$A_r$ implies that $A_r\cup\{x\}\notin K_r$. This is a
contradiction since the partition $\mathcal{Z} = \langle
Z_1,\dots, Z_r \rangle$ where $Z_\nu := X_\nu\setminus\{x\}, Z_r:=
A_r\cup\{x\}$ and $Z_j = X_j$ for $j\notin\{\nu, r\}$ clearly
violates the condition that $\mathcal{K}$ is collective
$r$-unavoidable.

\medskip
Let $V\subset [n]$. We say that a simplicial complex $K\subset
2^{[n]}$ is $V$-homogeneous if $S\in K \Leftrightarrow \phi(S)\in
K$ for each permutation $\phi : V\rightarrow V$ and each $S\subset
V$.

\medskip\noindent
{\it Claim~1.}  Each of the complexes $\{K_j\}_{j=1}^{r-1}$ is
$X$-homogeneous where $X = \bigcup_{j=1}^{r-1}~X_j  =[n] \setminus
A_r$.

\medskip\noindent {\it Proof of the Claim~1:} Let us show for
illustration that $K_1$ is $X$-homogeneous. This is deduced from
the observation that for each bijection $\phi : X\rightarrow X$,

\begin{equation}\label{eqn:a-b}
(a)\quad \phi(X_1)\notin K_1  \quad \mbox{\rm and }\quad (b)\quad
\phi(\partial(X_1))\subset K_1.
\end{equation}

This is obvious if $\phi(X_1) = X_1$. Moreover, it is sufficient
to establish (\ref{eqn:a-b}) in the case when $\phi$ is a
transposition, say $\phi(x_1) = x_2$ where $x_1\in X_1$ and
$x_2\in X_2$ (the case $x_2\in X_j$ for $j>2$ is treated
similarly).

\medskip\noindent
$(a)$\,  is equivalent to $(X_1\setminus\{x_1\})\cup\{x_2\}\notin
K_1 $. This is true since otherwise,
$$(X_1\setminus\{x_1\})\cup\{x_2\}\in K_1,\, X_2\setminus\{x_2\}\in
K_2,\, A_3\in K_3,\dots, A_{r-1}\in K_{r-1},\, A_{r}\in K_r $$
would be a disjoint family of sets covering all but $(r-2)$
elements of $[n]$ (contradicting (2) in Definition~3.4).

\smallskip\noindent
In order to prove $(b)$ let $X_1\setminus \{y\}$ be a facet of
$\partial(X_1)$ (the interesting case is $y\neq x_1$). Then, $X_1'
:= \phi(X_1\setminus \{y\}) = (X_1\setminus \{x_1, y\})\cup
\{x_2\}\in K_1$. Otherwise the disjoint collection,
$$X_1'\notin K_1,\, (X_2\setminus\{x_2\})\cup\{x_1\}\notin
K_2,\, X_3\notin K_3,\dots, X_{r-1}\notin K_{r-1},\,
X_r\cup\{y\}\notin K_r
$$
would violate the collective $r$-unavoidability of $\mathcal{K}$.
(Note that $(X_2\setminus\{x_2\})\cup\{x_1\} = \phi(X_2)\notin
K_2$ follows from $(a)$.) \hfill $\square$

\bigskip
Summarizing, we have so far established that for each $j=1,\dots,
r-1$ the restriction of $K_j$ on $X$ is the complex $X \choose
\leq m_j$ where $m_i$ is the cardinality of the set $A_i$. In
particular the sets $A_1, A_2,\dots, A_{r-1}$ can be replaced by
any disjoint family $A_1', A_2',\dots, A_{r-1}'$ of subsets of $X$
such that $\vert A_i'\vert = m_i$ for each $i$.

\bigskip\noindent
{\it Claim~2.} If $x\in X$ then $X_r\cup\{x\} = A_r\cup\{x\}\notin
K_r$.

\medskip\noindent
{\it Proof of the Claim~2:} By Claim~1 we can assume that $x\in
X\setminus\bigcup_{i=1}^{r-1}~A_i$. Then the assumption $A_r' :=
A_r\cup \{x\}\in K_r$ would contradicts the fact that $\mathcal{A}
= (A_1,\dots, A_r)$ is a maximal disjoint $r$-tuple of sets in
$[n]$ such that $A_i\in K_i$ for each $i=1,\dots, r$.  \hfill
$\square$

\medskip It follows from Claim~2 that either $X_r\cup\{x\} = A_r\cup\{x\}$ is a
$K_r$-blocker (this corresponds to the case (2) of the theorem) or
there exists a proper subset $S\subset X_r$ such that $S\cup\{x\}$
is a $K_r$-blocker. The following claim makes this observation
more precise by showing (eventually) that $S\subset X_r$ is unique
with this property (and in particular independent of $x$).

\bigskip\noindent {\it Claim~3.} Choose $x\in X$.  Let $X_r = S\uplus
C$ be a partition of $X_r$ such that $\{x\}\cup S$ is a
$K_r$-blocker, i.e.\ such that $\{x\}\cup S\notin K_r$ and
$\partial(\{x\}\cup S)\subset K_r$. Then $T\cup C\in K_i$ for each
$T\subset X\setminus\{x\}$ of cardinality $m_i$ where $i=1,\ldots,
r-1$. Moreover, $T\cup C$ is a facet (maximal simplex) of $K_i$.

\medskip\noindent
{\it Remark.} The case $C=\emptyset$ is NOT ruled out. As it will
turn out from the proof if $X_r = S'\uplus C'$ is another
decomposition such that $\{x\}\cup S'$ is a $K_r$-blocker then
$S'=S$ and $C'=C$.

\medskip\noindent
{\it Proof of the Claim~3:} Assume that $i=1$ (the proof in other
cases is analogous). Since the sets $A_1,\dots, A_{r-1}$ can be
replaced by any disjoint family $A_1', \dots, A_{r-1}'$ of subsets
of $X$ such that $\vert A_i'\vert = m_i$ for each $i$, we assume
that $T=A_1$. For a similar reason we can assume that $x\notin
\bigcup_{i=1}^{r-1}~A_i$.

Then $A_1\cup C\in K_1$ since otherwise,
$$
 A_1\cup C\notin K_1, \, X_2\notin K_2, \dots, X_{r-1}\notin
 K_{r-1},\, \{x\}\cup S\notin K_r
$$
would violate the collective $r$-unavoidability of $\mathcal{K}$.

\smallskip
Suppose that $A_1\cup C$ is not a facet of $K_1$. It follows that
$A_1\cup C\cup\{z\}\in K_1$ for some $z\in S$, hence
$(S\cup\{x\})\setminus\{z\}\in K_r$. This is a contradiction since
the family,
$$A_1\cup C\cup\{z\}\in K_1,\, A_2 \in
K_2,\, A_3\in K_3,\dots, A_{r-1}\in K_{r-1},\,
(S\cup\{x\})\setminus\{z\}\in K_r$$ is a disjoint family of sets
covering all but $(r-2)$ elements of $[n]$ (contradicting (2) in
Definition~3.4).  \hfill $\square$

\bigskip
It follows from Claim~3 that $K_1$ and $K_r$ can interchange
roles. More explicitly $B_r :=  S$ can be included in a disjoint
family $\{B_j\}_{j=1}^{r}$ (replacing the family
$\{A_i\}_{i=1}^{r}$) where $B_j = A_j$ for each $j=2,\dots, r-1$
and $B_1:= A_1\cup C$.

\smallskip
In light of Claim~1 each of the complexes $\{K_j\}_{j=2}^{r}$ is
$Y$-homogeneous where $Y = [n]\setminus (\{x\}\cup B_1)$ and $x$
is an arbitrary element in $[n]\setminus B_1$. Moreover the
decomposition $B_1 = A_1\uplus C$ corresponds to the decomposition
$X_r = A_r = S\uplus C$ in Claim~3.

\medskip
From here it is not difficult to conclude that for each
$i=1,\dots, r$ there is a decomposition $K_i \cong W_i\ast F$
where $W_i\cong {[n]\choose \leqslant m_i}$  and  $F$ is either
empty or $F = \Delta(C)$ is the simplex spanned by a finite,
non-empty set $C$. \hfill $\square$

\begin{cor}\label{cor:Bier-cor}  If $K =  K_1\ast_\Delta\dots\ast_\Delta K_r$ is a
Bier complex (Definition~\ref{dfn:r-Bier}) then either,

\begin{enumerate}
 \item[(1)]\hspace{0.5mm} $r=2$ and $K = K_1\ast_\Delta K_2 = K\ast_\Delta K^\circ$ is a Bier sphere, or
 \item[(2)]\hspace{0.5mm} $r\geq 3$ and $K = \Delta_{n,r}^{m_1,\dots, m_r; \mathbf{1}}$
is an optimal chessboard complex where $n = m_1+\dots +m_r + r-1$
(Section~\ref{sec:multi-chess}), or
 \item[(3)] \hspace{0.5mm} $r\geq 3$ and $K = \Delta \ast [r]^{\ast k}$ where $\Delta =\Delta_{n,r}^{m_1,\dots, m_r; \mathbf{1}}$
is an optimal chessboard complex and $[r]^{\ast k} =
[r]\ast\dots\ast [r]$ is the join of $k\geq 1$ copies of the
$0$-dimensional complex $[r]$.

 \end{enumerate}
{\rm (Note that {\rm (2)} is a formal `consequence' of {\rm (3)}
if we allow $k=0$.)}
\end{cor}

\textit{Proof.} Assume $r\geq 3$. It follows from
Theorem~\ref{thm:class} and Lemma~\ref{lemma:commute} that $K
\cong  E\ast F$ where (Section~\ref{sec:multi-chess})
$E=\Delta_{n,r}^{m_1,\dots, m_r; \mathbf{1}}$ is an optimal
chessboard complex and $F$ is either empty (the case (2)) or  $F =
(\Delta(C))^{\ast r}_\Delta$ for a non-empty set $C$ of
cardinality $\vert C\vert = k$. The proof is completed by the
observation that,
$$(\Delta(C))^{\ast r}_\Delta \cong (\{p\}^{\ast k})^{\ast
r}_\Delta \cong (\{p\}^{\ast r}_\Delta)^{\ast k}\cong [r]^{\ast
k}.$$

\section{Two perfect discrete Morse functions on the Bier sphere}
\label{sec:two-perfect}

\bigskip We illustrate the method of constructing DMF on deleted joins (by the method
of `migrating elements') first in the case of classical Bier
spheres.

It is known that Bier spheres are always shellable, see
\cite{BornerZiegler}. A method of Chari \cite{Chari} can be used
to turn this shelling into a perfect DMF on a Bier sphere. The
construction of the  `first perfect DMF' on a Bier sphere
(Section~\ref{sec:first-DMF}) essentially follows this path.

The `second perfect DMF' (Section~\ref{Sec SecondBierMorse})
differs from the first DMF, although the `migration rules' look
very similar. The advantage of the second DMF is that it can be
generalized to Alexander $r$-tuples and the associated
Bier complexes.

\subsection{First perfect DMF}\label{sec:first-DMF}

\par
We construct a discrete vector field $D_1$ on the Bier sphere
$Bier(K)$ in two steps:
\begin{enumerate}
    \item[(1)] We match  the simplices
$$\alpha=(A_1,A_2;B\cup i)  \hbox{ and } \beta=(A_1,A_2\cup i;B )$$ iff the following holds:
\begin{enumerate}
    \item $i<B,\ i<A_2$  \newline(that is, $i$ is smaller than all the entries of $B$ and $A_2$).
    \item $A_2\cup i \in K^{\circ}$.
\end{enumerate}
\end{enumerate}

\medskip\noindent
Before we pass to step 2, let us observe that the non-matched
simplices are labelled  by $(A_1,A_2;B\cup i)$ such that $A_2\in
K^{\circ}$, but $A_2\cup i \notin K^{\circ}$. As a consequence,
for non-matched simplices $A_1\cup B\in K$.

\smallskip
\begin{enumerate}
    \item[(2)] In the second step we match together the simplices
$$\alpha=(A_1,A_2;B\cup j)  \hbox{ and } \beta=(A_1\cup j,A_2;B)$$ iff the following holds:
\begin{enumerate}
\item None of the simplices $\alpha$ and $\beta$ is matched in the
first step.
    \item $j>B,\ j>A_1$.
    \item $A_1\cup j \in K$.
\end{enumerate}
\end{enumerate}
\noindent Observe that the condition (c)  always holds (provided
that the condition (a) is satisfied), except for
    the case  $B=\emptyset$.

\begin{lemma} The discrete vector field $D_1$ is a discrete Morse function on the Bier sphere $Bier(K)$.
\end{lemma}

\textit{Proof.} Since $D_1$ is (by construction) a discrete vector
field,  it remains to check that there are no closed gradient
paths. Observe that in each pair of simplices in the discrete
vector field $D_1$ there is exactly one \textit{migrating
element}. More precisely, in the case (1) the element $i$
{migrates} to $A_2$, and in  the case (2) the element $j$
{migrates} to $A_1$.

The lemma follows from the observation that (along a gradient
path) the values of the migrating element that move to $A_2$
strictly decreases. Similarly, the values of migrating elements
that move to $A_1$ can only increase.

\medskip
Let us illustrate this observation by an example. Assume we have a
fragment of a gradient path that contains two matchings of type 1.
We have:

$$(A_1\cup k ,A_2;B\cup i ) \rightarrow (A_1\cup k ,A_2\cup i;B)\rightarrow $$$$(A_1,A_2\cup i;B\cup k ) \rightarrow (A_1,A_2\cup k \cup i;B)$$
The migrating elements here are $i$ and $k$. The definition of the
matching $D_1$ implies $k < i$. Otherwise $(A_1,A_2\cup i;B\cup k
)$ is matched with $(A_1,A_2;B\cup k\cup i )$, and the path would
terminate after its second term. \qed

\medskip

It is not difficult to see that there are precisely two critical
simplices in $D_1$:

\begin{enumerate}
    \item An $(n-2)$-dimensional simplex,  $$(A_1,A_2; i )$$
where $A_1<i<A_2$,  (this condition describes this simplex
uniquely, in light of the fact that $A_1\in K$ and $A_2\in
K^{\circ}$),
    \item and the $0$-dimensional simplex, $$(\emptyset, \{1\} ;\{2,3,4,...,n\}).$$
\end{enumerate}
\noindent (Here we make a simplifying assumption that $\{1\}\in
K^\circ$, which can be always achieved by a re-enumeration, except
in the trivial case $K^\circ =\{\emptyset\}$.)

\subsection{Second perfect DMF}\label{Sec SecondBierMorse}
The construction of the second discrete vector field $D_2$ is also
in two steps:

\smallskip
The first step remains the same:
\begin{enumerate}
    \item We match  the simplices
$$\alpha=(A_1,A_2;B\cup i)  \hbox{ and } \beta=(A_1,A_2\cup i;B )$$ iff the following holds:
\begin{enumerate}
    \item $i<B,\ i<A_2$  \newline(that is, $i$ is smaller than all elements in $B$ and $A_2$).
    \item $A_2\cup i \in K^{\circ}$.
\end{enumerate}

\bigskip

Before we pass to the second step, let us remind ourselves that
the non-matched simplices are labelled  by $(A_1,A_2;B\cup i)$ such
that $A_2\in K^{\circ}$, but $A_2\cup i \notin K^{\circ}$. As a
consequence, for non-matched simplices $A_1\cup B\in K$.

\bigskip

    \item We match together the simplices
$$\alpha=(A_1,A_2,B\cup i \cup j )  \hbox{ and } \beta=(A_1\cup j ,A_2,B\cup  i )$$ iff the following holds:
\begin{enumerate}
\item None of the simplices $\alpha$ and $\beta$ was matched in the
first step, i.e.
  \newline  $i<j,\  i<B,\  i<A_2$, and $i\cup A_2 \notin K^{\circ}$.
    \item $j<B,\ j<A_1\setminus [1,i]$.
    \item $A_1\cup j \in K^{\circ}$.
\end{enumerate}
\end{enumerate}

\noindent Note that the condition (c) is always satisfied
(provided that the condition (a) above holds).

\medskip

We omit the proof that $D_2$ is indeed a discrete Morse function
since a more general fact will be established in the proof of
Theorem~\ref{ThmMain} (Section~\ref{SecUnavoid}).

\medskip

Finally we observe that, with the same simplifying assumption
$\{1\}\in K^\circ$, the discrete vector fields $D_2$ and $D_1$
have the same critical simplices:
\begin{enumerate}
    \item $(A_1,A_2, i )$
such that $A_1<i<A_2$
    \item and $(\emptyset, \{1\}; \{2,3,4,...,n\}).$

\end{enumerate}

\section{ Proof of Theorem \ref{ThmMain}}\label{SecUnavoid}

The  proof of Theorem~\ref{ThmMain} is based on the construction
of a discrete Morse function $D$ on the deleted join
$\mathcal{K}^{\ast r}_\Delta$, where $\mathcal{K} = \langle
K_1,...,K_r\rangle$ is a collective $r$-unavoidable collection of
complexes.

\medskip We will demonstrate that:
 \begin{itemize}
   \item If the $r$-tuple $\mathcal{K}$ is  collective $r$-unavoidable, then the critical simplices of the discrete Morse field $D$
    may appear only starting with dimension $n-r$
    (except for the unique $0$-dimensional simplex).
       This observation immediately implies the connectivity bound in Theorem~\ref{ThmMain}, part (1).
   \item Under the stronger hypothesis that $\mathcal{K}$ is an Alexander
   $r$-tuple, the discrete vector field $D$ has a single
   $0$-dimensional critical simplex, while all other critical simplices have one and the same dimension $n-r$.
Theorem \ref{ThmMain} (part (2)) is an immediate consequence.
Moreover, a direct dimension count will establish the purity of
the complex $\mathcal{K}^{\ast r}_\Delta$.
 \end{itemize}

\medskip
As in Section~\ref{sec:collectively}, a simplex $\beta =
A_1\uplus\dots\uplus A_r \in \mathcal{K}^{\ast r}_\Delta$ is in
the `partition notation' recorded as $\beta = (A_1,\dots, A_r; B)$
where $B = [n]\setminus \cup_{i=1}^r~A_i$. More explicitly, an
ordered partition $(A_1,A_2,...,A_r;B)$  of $[n]$ into $r+1$
parts, corresponds to a simplex in $\mathcal{K}^{\ast r}_\Delta$
if and only if,

\begin{enumerate}
 \item[(1)] $A_i  \in K_i$ for each $i=1,\dots, r$;
 \item[(2)] $\cup\{A_i\}_{i=1}^r\neq\emptyset$, meaning that the partition
 $(\emptyset, \dots, \emptyset , [n])$ is excluded.
\end{enumerate}

\medskip Observe that the dimension of a simplex $\beta = (A_1,\dots, A_r; B)$
is determined by the cardinality of $B$, indeed   ${\rm
dim}(\beta)=n-|B|-1.$

Moreover, a facet of a simplex $\beta = (A_1,A_2,...,A_r;B)$ is
obtained by moving (we also say `migrating') an element from one
of the sets $A_i$ to $B$.
  For example, $(\{1,2\},\{6\},\{5\};\{3,4, 7\})$ is a facet of
  $(\{1,2\},\{6,7\},\{5\};\{3,4\})$ obtained by the migration of the element $7\in A_2$.

\bigskip

\subsection*{Construction of the discrete Morse function $D$} The
discrete vector field $D$ is described by a step-by-step
construction, generalizing the construction of the discrete vector
field $D_2$ from Section~\ref{Sec SecondBierMorse}.

\medskip  In the first step we match the simplices,
$$\alpha=(A_1,A_2,...,A_r;B\cup i_{1} )  \hbox{ and } \beta=(A_1\cup i_{1},A_2,...,A_r;B)$$ iff the following holds:
\begin{enumerate}
    \item $i_{1}<B,\ i_{1}<A_1$;
    \item $A_1\cup i_{1} \in K_1$.
\end{enumerate}
In other words a simplex $\alpha = (A_1,A_2,...,A_r;B')$ is
matched (if possible) with the simplex $\beta =
(A_1',A_2,...,A_r;B)$ obtained from $\alpha$ by migrating the
minimum $i_1$ of the set $B' = B\cup i$ into $A_1' = A_1\cup i$
(provided $i_1< A_1$ and $A_1'\in K_1$).

\medskip Observe that many simplices are matched already in this
step. Indeed, for $\alpha = (A_1',A_2,...,A_r;B')$ let $i_1 = {\rm
min}(A_1'\cup B')$. If $i_1\in B'$ then $\alpha$ is matched with
the simplex obtained by migrating $i_1$ from $B'$ to $A_1'$. If
$i_1\in A_1'$, $\alpha$ is obtained by the migration of $i_1$ from
its facet $\gamma = (A_1,A_2,...,A_r;B)$, where $A_1 =
A_1'\setminus \ i_1$ and $B = B'\cup i_1$.

\medskip
The remaining non-matched simplices $(A_1,A_2,...A_r;B\cup i_1) $
fall into two types:

\begin{enumerate}
  \item The first type:
\begin{enumerate}

  \item[] $i_1<B$, $i_1<A_1$ and $A_1\cup i_1 \notin K_1$.
\end{enumerate}

  \item
The second type:

\begin{enumerate}
  \item[] $B=\emptyset$ and $A_1 =\emptyset$.
\end{enumerate}
\end{enumerate}

\noindent Here we declare that the non-matched simplices {\bf of
the second type} {\em will not participate} in matching in later
steps of the construction, i.e.\ they will contribute to the
critical simplices of $D$.

\medskip
There is a single $0$-dimensional non-matched simplex, $(\{1\},
\emptyset,..., \emptyset ; \{2,\dots, n\})$. Here (as in
Section~\ref{sec:two-perfect}) we make a simplifying
(non-essential) assumption that $\{1\}\in K_1$. (This condition
can be easily satisfied by choosing a different linear order on
$\mathcal{K}$ and $[n]$, if necessary.)

\medskip
We continue the construction by trying to migrate elements from
$B$ into $A_2$ (in the second step), into $A_3$ (in the third
step), etc. Assume, inductively, that the first $(k-1)$-steps of
the construction are completed.

\bigskip In the $k$-th step of the construction we match the
simplices,
$$\alpha=(...,A_k,...,A_r; B\cup i_1... \cup i_{k-1}\cup i_k)  \hbox{ and } \beta=(...,A_k\cup i_k,...,A_r;B\cup i_1... \cup i_{k-1})$$
iff the following holds:
\begin{enumerate}
\item $\alpha$ and $\beta$ are non-matched simplices of {\bf the
first type} in all preceding steps,

    \item $i_k<B,\ i_k< A_k\setminus [1,i_{k-1}]$,
    \item $A_k \cup i_k \in K_k$.
\end{enumerate}

The remaining non-matched simplices   $(A_1,A_2,...,A_k,...A_r;B\cup i_1... \cup i_{k-1}\cup i_k) $  again fall into two types:
\begin{enumerate}
  \item The first type:
\begin{enumerate}
  \item the simplex is a first type non-matched simplex on steps $1,...,k-1$,
  \item $i_k<B$, $i_k<A_k \setminus [1,i_{k-1}]$,
  \item $A_k\cup i_k \notin K_k$.
\end{enumerate}

  \item
The second type:

\begin{enumerate}
  \item the simplex is a first type non-matched simplex on steps $1,...,k-1$,
  \item $B=\emptyset$, $A_k \subset [1,i_{k-1}]$.
\end{enumerate}
\end{enumerate}

\noindent (As before we declare that the non-matched simplices of
the second type never participate in subsequent matchings.)

\medskip
From the assumption that $\mathcal{K}$ is a collective
$r$-unavoidable collection of complexes we conclude that on the
$r$-th step there are no non-matched simplices of the first type.
From here we deduce that the cardinality of $B$ for critical
simplices can vary from $0$ to $r-1$, and in particular the
dimension of any  critical simplex is greater or equal than $n-r$.
(The only exception being of course the $0$-dimensional critical
simplex $(\{1\},\emptyset,...,\emptyset;\{2,...,n\})$.)

%\textbf{If we relax the non-avoidability condition we conclude that cardinality of $B$ for critical
%simplices varies from $r-1-h$ to $r$. Therefore their dimension varies from $n-r+1$ to $n-r+h$.}

\bigskip

\subsection*{An alternative description of the DMF}
It may be instructive to summarize the construction of the
discrete Morse function $D$ in the form of an `algorithm' which
describes the matching and lists the critical simplices.

For this purpose we introduce an operator $\textbf{a}$ which takes
simplices $$(A_{1},\dots,A_{r};B) \in
K_{1}*_{\Delta}\dots*_{\Delta}K_{r}$$ and maps them  to strictly
increasing $r$-tuples, $$\textbf{a}=(a_{1}< a_{2}< \dots<
a_{r})\in (\mathbb{N} \cup \{\infty\})^r,$$ by the following rule:

\begin{description}
    \item[1] \quad $a_{1}:=\min (B\cup A_{1})$; \,
if $B\cup A_{1}=\emptyset$ then $a_1=...=a_{r}:=\infty$.
    \item[2] \quad $a_{2}:=\min ((B\cup A_{2})\backslash
    [1,a_{1}])$;

\qquad \qquad if $(B\cup A_{2})\backslash [1,a_{1}] = \emptyset$
then   $a_k:=\infty$ for all $k\geq 2$.$$\dots$$

   \item[i] \quad
$a_{i}:=\min ((B\cup A_{i})\backslash [1,a_{i-1}]);$

\qquad \qquad if $((B\cup A_{i})\backslash [1,a_{i-1}])=\emptyset$
then   $a_k:=\infty$ for all $k\geq i$. $$\dots$$
    \item[r] \quad
$a_{r}:=\min \{(B\cup A_{r})\backslash [1,a_{r-1}]\};$

\qquad \qquad  if $((B\cup A_{r})\backslash
[1,a_{r-1}])=\emptyset$ then   $a_r:=\infty$.
\end{description}

\bigskip

%\begin{dfn}
\noindent  We say that an element $a_{j}$ of the $r$-tuple
$\textbf{a}(A_{1},\dots,A_{r};B)$ is \textit{potentially movable} if $a_{j}\neq \infty$. A potentially movable element $a_j$ is  \textit{movable} if:

\begin{enumerate}
\item either $a_{j} \in B$ and $A_{j}\cup a_{j} \in K_{j}$,
\item or $a_{j} \in A_j$.
\end{enumerate}
The \textbf{standard move} of a movable element $a_{j}$ is the
matching of:
\begin{enumerate}
  \item either $(A_{1},\dots,A_{j},\dots,A_{r};B) \rightarrow (A_{1},\dots,A_{j} \cup a_{j},\dots,A_{r};B\backslash a_{j}),$
  \item or $(A_{1},\dots,A_{j}\setminus a_j ,\dots,A_{r};B\cup a_{j}) \rightarrow (A_{1},\dots,A_{j},\dots,A_{r};B).$
\end{enumerate}

\medskip

The following procedure finds the corresponding pair (if any) for
each simplex. If the simplex is not matched, the algorithm reports
that it is critical.

\medskip

\subsection*{Matching Algorithm} A simplex
$(A_{1},\dots,A_{r};B)\in K_{1}*_{\Delta}\dots*_{\Delta}K_{r}$ is
matched with the simplex obtained by the standard move of  the
\textit{minimal movable element} in
$\textbf{a}(A_{1},\dots,A_{r};B)$. If there are no movable
elements, the simplex is critical.

\medskip
\begin{prop} The "Matching Algorithm" describes a discrete Morse function $D$.
\end{prop}
It is clear that $D$ is a discrete vector field. The proof of the
acyclicity follows from the following lemmas.

\begin{lemma}\label{lemma:decreases}
 Assume the lexicographic order on the set $(\mathbb{N} \cup
\{\infty\})^r$. Then the function $\textbf{a}$ decreases
(non-strictly) along any gradient path of the discrete vector
field $D$ described by the "Matching Algorithm".
\end{lemma}

\textit{Proof.} For any gradient path,
$$\alpha_0^p\rightarrow \beta_0^{p+1}\rightarrow \alpha_1^p\rightarrow \beta_1^{p+1}\rightarrow ..., $$
we observe that $\textbf{a}(\alpha_i^p)=\textbf{a}(
\beta_i^{p+1})$  and $\textbf{a}( \beta_i^{p+1})\leq
\textbf{a}(\alpha_{i+1}^p)$.\qed

\medskip
It immediately follows from Lemma~\ref{lemma:decreases} that the function $\textbf{a}$ must be constant along a cyclic gradient path (if it exists). In particular, along such a path the set of all potentially movable elements remains the same. The following lemma rules out this possibility.

\begin{lemma}\label{lemma:acyclic}
If the function $\textbf{a} = (a_{1}, a_{2}, \dots, a_{r})$ is
constant along a gradient path, then the path is acyclic.
\end{lemma}

\textit{Proof.} Let us inspect a typical fragment of the gradient
path,  
\begin{equation}\label{eqn:fragment}
 \alpha_0^p\longrightarrow \beta_0^{p+1}\dashrightarrow \alpha_1^p,
\end{equation}
which is more explicitly recorded as the path,
$$(...,A_r;B\cup a_k) \rightarrow (...,A_k\cup  a_k,... ;B)\dashrightarrow
(...,A_k\cup a_k,..., A_m\setminus \nu,...;B\cup \nu).$$
 If $\nu<a_k$, the value of $a_k$ would change, contrary to the assumption that $\textbf{a}$ is constant along the path. In the case $\nu> a_k$ there are two possibilities.

The first possibility is that (by the matching algorithm) there is a matching,
\[
(A_1\,...,A_k,..., A_m\setminus \nu,...,A_r; B\cup a_k\cup \nu)
\longrightarrow
(A_1\,...,A_k\cup a_k,..., A_m\setminus \nu,...;B\cup \nu)
 \]
 or in other words a matching $\gamma^{p-1} \longrightarrow \alpha_1^p$,
 which would guarantee that the gradient path (\ref{eqn:fragment}) terminates at $\alpha_1^p$.

This happens precisely if $a_k$ is the minimal movable element in,
\begin{equation}\label{eqn:second-case}
\textbf{a}(A_1\,...,A_k\cup \{a_k\},..., A_m\setminus \nu,...,A_r;B\cup \nu).
\end{equation}
The only possible scenario when $a_k$ is not the minimal movable element in (\ref{eqn:second-case}) is when $m<k$ and the element $a_m$ happens to be {\it  movable} (as a consequence of $\{a_m\}\cup(A_m\setminus\{j\})\in K_m$). 
%This is indeed possible as shown by the Example~\ref{Ex:gradient-path}.

\medskip
Summarizing we observe that in the `worst case scenario' the minimal movable element $a_m$ of the simplex $\alpha_1^p$ in (\ref{eqn:fragment}) is strictly smaller than the minimal movable element $a_k$ of the simplex $\alpha_0^p$. It follows that if this case persists, then the minimal movable element decreases along the path and the path must be acyclic.  \qed

\medskip

To establish the second statement in Theorem~\ref{ThmMain}, we
need the second half of Theorem \ref{ThmWedge}. By
Proposition~\ref{prop:Pure} if $\mathcal{K}$ is an Alexander
$r$-tuple, then the complex $\mathcal{K}^{\ast r}_\Delta$ is pure
$(n-r)$-dimensional.

\medskip
We end the proof with an efficient, combinatorial description of
critical cells of the discrete Morse function $D$. The
Corollary~\ref{cor:korolar} is a consequence of the well known
fact that the spheres in the wedge decomposition of
$Bier(\mathcal{K})$ are in one-to-ne correspondence with the
critical cells of $D$.

\begin{enumerate}
    \item[1.]  An $(n-r)$-dimensional simplex $$(A_1,A_2,\dots,A_r;i_1\cup i_2 \cup\dots\cup i_{r-1}) \ \ \hbox{  with  } \ \ i_1<i_2<...<i_{r-1}$$
is a critical simplex of the discrete Morse function $D$ if and
only if:
\begin{enumerate}
\item $A_1$ avoids the segment $[1,i_1]$,
\item $A_k$ avoids the segment $[i_{k-1},i_k]$  for $k\in [2,r-1]$,
\item $A_{r-1}$ avoids the segment $[i_{r-2},i_{r-1}]$,
\item $A_{k}\cup i_{k} \notin K_{k}$ for $k\in [1,r-1]$,
\item $A_{r}$ avoid the segment $[i_{r-1}, n]$.
\end{enumerate}
\bigskip
    \item[2.] There is a single $0$-dimensional simplex:
    $$( \{1\},\emptyset,\dots,\emptyset; \{2,3,4,...,n\}).$$

\end{enumerate}
With this observation we complete the proof of Theorem
\ref{ThmMain} (Corollary~\ref{cor:korolar}).

\begin{Ex}{\em  The chessboard complex
$\Delta_{5,3}^{1,1,1;\bf{1}}$ is the  Bier {complex}
associated to the Alexander $3$-tuple $\mathcal{K} = \langle K_1,
K_2, K_3 \rangle$, where $K_{1}=K_{2}=K_{3}$ is the
$0$-dimensional skeleton of the $4$-dimensional simplex
$\Delta([5])$. Then the critical simplices of the discrete Morse
function $D$ constructed in the proof of Theorem~\ref{ThmMain} are
the following:

\medskip\noindent
$(A_{1}, A_{2}, A_{3}; B)=(4,5,2;\{1,3\}) ,$ $(5,4,2;\{1,3\}),$
$(2,5,3;\{1,4\}),$ $(3,5,3;\{1,4\}),$ $(4,5,1;\{2,3\}),$
$(5,4,1;\{2,3\}),$ $(3,5,1;\{2,4\}),$ $(5,1,3;\{2,4\}),$
$(3,1,4;\{2,5\}),$ $(4,1,3;\{2,5\}),$ $(5,2,1;\{3,4\}),$
$(5,1,2;\{3,4\}),$ $(4,2,1;\{3,5\})$ $(4,1,2;\{3,5\})$, and
$(1,\emptyset,\emptyset;\{2,3,4,5\})$. }
\end{Ex}

%{\color{blue}
%\begin{Ex}\label{Ex:gradient-path}
%The following gradient path illustrates a possible behaviour of the function %$\textbf{a}$ (illustrating Lemma~\ref{lemma:acyclic}) in the case of the complex %$\Delta_{5,3}^{1,1,1;\bf{1}}$.
%\begin{eqnarray*}
%(5,\emptyset, 4; \{1,2,3\}) \longrightarrow (5,2,4; \{1,3\}) \dashrightarrow %(\emptyset, 2, 4; \{1,3,5\}) \longrightarrow \\
%\longrightarrow  (1,2,4;\{3,5\}) \dashrightarrow (1,2,\emptyset; \{3,4,5\}) %\longleftarrow (\emptyset, 2, \emptyset; \{1,3,4,5\})
%\end{eqnarray*}
%\end{Ex}
%}

\section{Enumeration of critical simplices\\ for optimal chessboard
complexes}\label{SecNumb}

Optimal  chessboard complexes (see
Section~\ref{sec:multi-chess} and Example~\ref{ex:skeletons} in
Section~\ref{sec:collectively}) are our key examples of Alexander
$r$-tuples for $r\geq 3$. In this section we enumerate critical
simplices of $\Delta^{m_1,\ldots,m_r;1}_{n,r}$ for $n =
m_1+\cdots+m_r+r-1$.

\medskip
For a given simplex $\beta\in \Delta^{m_1,\ldots,m_r;1}_{n,r}$,
let us encode the set of free columns $1\leqslant
x_1<\ldots<x_{r-1}\leqslant n$ as an $(r-1)$-tuple
$\mathbf{x}=(x_1,x_2,\ldots,x_{r-1})$. Observe that if (in the
notation of Section~\ref{SecUnavoid}) $\beta = (A_1,\dots, A_n;B)$
then $B = \{x_i\}_{i=1}^{r-1}$.

\medskip
Let $\mathbf{b}=(b_1,b_2,\ldots,b_r)$ denote the sequence that
counts the number of rooks (in all rows) between consecutive free
columns, i.e.
$$b_1=x_1-1, b_2=x_2-x_1-1, \ldots, b_r=r-x_r.$$

Let $b_{i,j}$ denote the number of rooks in the $j^{th}$ column
between columns $x_{i-1}$ and $x_i$ (with the obvious
interpretation of numbers $b_{i1}$ and $b_{ir}$). We know from the
"critical simplices criterion"  (found at the end of the proof of
Theorem~\ref{ThmMain} in Section~\ref{SecUnavoid}), that
$b_{ii}=0$ for all $i$. For a given $\mathbf{x}$, all possible
numbers of rooks between free columns in critical simplices (we
ignore for a moment the order or rooks),
corresponds to all non-negative $r\times r$ matrices, $$B=\left(%
\begin{array}{cccc}
b_{11} & b_{12} & \cdots & b_{1r}  \\
b_{21} & b_{22} & \cdots &  b_{2r} \\
\cdots &\cdots & \cdots & \cdots \\
b_{r1} & b_{r2} & \cdots &  b_{rr}\\
\end{array}%
\right)$$ with non-negative integers such that:
$$b_{11}=\cdots=b_{rr}=0,\, B\cdot \mathbf{1}=(m_1,\ldots,m_r),\,
 \mathbf{1}^t\cdot B=(b_1,b_2,\ldots,b_r).$$
Denote the number of such matrices by $R_\mathbf{x}$. If the
number of rooks between two consecutive free columns for each row
is fixed, the number of all configurations is
$${b_1\choose b_{11},b_{21},\ldots, b_{r1}}{b_2\choose b_{12}, b_{22},\ldots, b_{r2}}\cdots
{b_r\choose b_{1r},b_{2r},\ldots, b_{rr}}.$$

Therefore, the number of all critical simplices of
$\Delta^{m_1,\ldots,m_r;1}_{n,r}$ is
$$\sum_{1\leqslant x_1<\ldots<x_{r-1}\leqslant r}R_\textbf{x}{b_1\choose b_{11},b_{21},\ldots, b_{r1}}{b_2\choose b_{12},b_{22},\ldots, b_{r2}}\cdots
{b_r\choose b_{1r},b_{2r},\ldots, b_{rr}}$$

\section{Discrete Morse function for a long chessboard complex}\label{secLong}

The `long' chessboard complexes (described in
Section~\ref{sec:multi-chess}) are not collective $r$-unavoidable
complexes, let alone Alexander $r$-tuples. However, the
construction of the discrete Morse function, described in
Section~\ref{SecUnavoid}, is sufficiently general and versatile to
be applied in this case as well. This is very interesting since
the existence of a perfect Morse function on this complex provides
an alternative proof of the (critical case) of Theorem~3.2 from
\cite{jvz-1}. Recall that this result paved the way for some new
Tverberg-Van Kampen-Flores type results, including the Theorem~1.2
from \cite{jvz-2}.

\medskip
Recall that a multiple chessboard complex
$K=\Delta_{n,r}^{m_{1},\dots,m_{r};\bf{1}}$ is `long' if $n >
m_1+\dots + m_r + r-1$.
Following essentially the construction of the matching described in
Section~\ref{SecUnavoid}, one obtains the discrete Morse function
$D$ on $K$ which has the following critical simplices.

With the exception of  the unique $0$-dimensional critical
simplex, all other critical simplices are described as the
configurations $(A_{1}, A_{2},\dots, A_{r};B)$ for which there
exist elements $i_{1} < i_{2} <\dots < i_{r} $ in $B$ such that
the following conditions are satisfied:

\begin{enumerate}
\item $A_1$ avoids the segment $[1,i_1]$,
\item $A_k$ avoids the segment $[i_{k-1},i_k]$  for $k\in [2,r]$,

\item $A_{k}\cup i_{k} \notin K_{k}$ for $k\in [1,r-1]$,
\item $B\backslash \{i_{1},\dots,i_{r}\}>i_{r}$.
\end{enumerate} The condition (3) implies that all the critical
simplices have one and the same dimension $(m_{1}+\dots+m_{r}-1)$.
%

%If we set $i_{0} := 0$ then the conditions can be reformulated like

%\textbf{For any $j\in [1,r]$ we have $A_{j} \cap [i_{j-1},i_{j}] = \emptyset$.}

\begin{Ex}{\em
The complex $\Delta_{4,2}^{1,1;\bf{1}}$ has $5$ critical simplices
of dimension $1$:

$(4,3,\{1,2\})$ $(3,4,\{1,2\})$ $(2,4,\{1,3\})$ $(4,1,\{3,2\})$
$(3,1,\{4,2\})$. }
\end{Ex}

The existence of a perfect discrete Morse function on the long,
multiple chessboard complex
$\Delta_{n,r}^{m_{1},\dots,m_{r};\bf{1}}$ provides an alternative
proof of the following theorem from \cite{jvz-1}. (Two other
proofs, both of them comparatively complex and non-trivial, relied
respectively on a shelling construction, and the Nerve Lemma.)

\begin{thm}{\em (\cite[Theorem~3.2]{jvz-1})}
The long chessboard complex is homotopy equivalent to a wedge of
$(m_{1}+\dots+m_{r}-1)$-dimensional spheres.
\end{thm}

\subsection{Enumeration of critical simplices for a long chessboard complex}
\label{sec:enumeration}

In this section we enumerate the critical simplices in the long
chessboard complex $\Delta^{m_1,\ldots,m_r;1}_{n,r}$.

\medskip We use the notation as in Section~\ref{SecNumb}. Recall
that (for all $i=1,2,\ldots,r$) the integer $b_{ij}$ evaluates the
number of rooks in the $i^{th}$ row between $(j-1)^{th}$ and
$j^{th}$ free column (for $j=1,2,\dots,r$). The distributions of
rooks between columns (we again ignore for a moment the exact
positions of the rooks) is encoded by the matrix $B\in
Mat_{r,r+1}(\mathbb{N}_0)$ where
$$ b_{11}=\cdots=b_{rr}=0,\, B\cdot \mathbf{1}=(m_1,\ldots,m_r).$$

Also, in this case we have $n-r-m_1-m_2-\cdots-m_r$ free columns,
and all of them are positioned behind the $r^{th}$ free column.

Simply by counting all partitions of the corresponding multisets,
we obtain the following formula for the number of all critical
simplices in $\Delta^{m_1,\ldots,m_r;1}_{n,r}$,
$$\sum_{\begin{array}{c}
  \scriptstyle{B\in
\mathbb{M}_{r,r+1}(\mathbb{N}_0),}\\
 \scriptstyle{ B\cdot \mathbf{1}=\mathbf{m}, b_{11}=\cdots=b_{rr}=0}\\
\end{array}}{n-r-\sum_{i}m_i+\sum_{j}b_{jr}\choose \sum_{j}b_{jr}}
\prod_{i=1}^{r}{b_{11}+b_{21}+\cdots+b_{r1}\choose b_{11},
b_{21},\ldots,b_{r1}}
$$
For example, if $r=2$ we have that the number of critical simplices
is
$$\sum_{b_{31}=1}^{m_1}\sum_{b_{32}=1}^{m_2}{b_{13}+b_{23}\choose b_{13}}
{n-2-m_1-m_2+b_{13}+b_{23}\choose b_{13}+b_{23}}.$$

%\section{Concluding remarks}
%\textcolor{red}{
%It seems plausible that one can compute homology groups for  a
%deleted product of $r$-tuples complexes that are "almost Alexander",
%that is, with critical simplices appearing  in two dimensions
%(rather than in just one) by counting gradient paths and describing
%Morse chain complex.}

%\textcolor{red}{Another interesting question is to develop a machinery that either  produces
%or characterizes  all existing Alexander $r$-tuples. Since Alexander $r$-tuples generalize simple games,
%one expects some technique similar to that from \cite{SimpleGames}.}

\end{document}